\documentclass[final]{l4dc2024}

\title[Probabilistic ODE Solvers for 
Numerical Optimal Control]
{Probabilistic ODE Solvers for Integration Error-Aware \\
Numerical Optimal Control}
\usepackage{times}
\usepackage{cleveref}
\usepackage{caption}
\usepackage{graphicx}

\usepackage{acronym}
\usepackage{doi} %

\usepackage{pgf}

\author{%
 \Name{Amon Lahr} \Email{amlahr@ethz.ch} \\
 \addr ETH Zürich
 \AND
 \Name{Filip Tronarp} \Email{filip.tronarp@matstat.lu.se} \\
 \addr Lund University
 \AND
 \Name{Nathanael Bosch} \Email{nathanael.bosch@uni-tuebingen.de} \\
 \addr University of Tübingen, Tübingen AI Center
 \AND
 \Name{Jonathan Schmidt} \Email{jonathan.schmidt@uni-tuebingen.de} \\
 \addr University of Tübingen, Tübingen AI Center
 \AND
 \Name{Philipp Hennig} \Email{philipp.hennig@uni-tuebingen.de} \\
 \addr University of Tübingen, Tübingen AI Center
 \AND
 \Name{Melanie N. Zeilinger} \Email{mzeilinger@ethz.ch} \\
 \addr ETH Zürich
  }

\begin{document}

\acrodef{MPC}{model predictive control}
\acrodef{ODE}{ordinary differential equation}
\acrodef{OCP}{optimal control problem}
\acrodef{IVP}{initial value problem}
\acrodef{MLE}{maximum likelihood estimate}
\acrodef{EKF}{Extended Kalman Filter}
\acrodef{EKS}{Extended Kalman Smoother}
\acrodef{(I)EKS}{(Iterated) Extended Kalman Smoother}
\acrodef{IEKS}{Iterated Extended Kalman Smoother}
\acrodef{RTS}{Rauch-Tung-Striebel Smoother}

\maketitle

\begin{abstract}%
  Appropriate time discretization is crucial for 
  real-time applications of numerical optimal control, such as nonlinear model predictive control.
  However, 
  if
  the discretization error strongly depends on the applied control input, meeting accuracy and sampling time requirements simultaneously can be challenging using classical discretization methods. 
  In particular, 
  neither fixed-grid nor adaptive-grid discretizations may be suitable,
  when they
  suffer from 
  large integration error 
  or 
  exceed the prescribed sampling time,
  respectively. 
  In this work, 
  we take a first step at 
  closing this gap 
  by 
  utilizing probabilistic 
  numerical integrators
  to approximate the solution of the 
  initial value problem,
  as well as the computational uncertainty associated with it,
  inside the 
  \ac{OCP}.
  By taking the viewpoint of probabilistic numerics
  and propagating
  the numerical uncertainty 
  in the cost, 
  the \ac{OCP} is reformulated such that the optimal input
  reduces the computational uncertainty insofar as it is beneficial for the control objective. 
  The proposed approach 
  is illustrated
  using a numerical example,
  and 
  potential benefits and limitations are discussed.
  \end{abstract}


\acresetall

\section{Introduction}

For
direct optimal control, 
numerical integration methods are a crucial component for converting a continuous-time system model of \acp{ODE} into accurate discrete-time predictions based on which the input signal can be optimized. 
While negligible discretization error can generally be achieved by choosing a sufficiently fine time discretization, 
it can be computationally infeasible
in the \ac{MPC} setting, %
where an \ac{OCP} is to be solved online, i.e., in receding horizon fashion and at a fixed sampling interval.  
This is complicated by the fact that the accumulated discretization error at each time step may depend on the chosen input signal, directly due to the control input affecting the vector field, 
or indirectly due to changing stiffness properties between different local solutions of the resulting state trajectory. 

For such challenging scenarios, 
standard grid discretization approaches may fail to 
achieve negligible integration error 
given a limited computational budget.
First, a 
uniform time discretization may lead to either unacceptable integration error for a coarse grid or 
excessive computational requirements for a fine grid.
Second, 
a fixed, non-uniform time discretization~\citep{qin1997overview,tondel_complexity_2002,cagienard_move_2007,shekhar_optimal_2015,yu_stable_2016,chen_efficient_2020} 
may reduce discretization error without adding computational overhead~\citep{quirynen_multiple_2015}. 
The offline design of fixed, non-uniform grids has been considered 
by \cite{lazutkin_approach_2018} and \cite{wolff_nonlinear_2022}.
\cite{lishkova_multirate_2022}
consider
multi-rate time discretization for dynamical systems with varying time scales within the MPC.
Still, the discretization accuracy of these approaches may vary as they do not consider the dependence of the integration error on the applied input. 
Third, 
adaptive mesh refinement methods 
can achieve a pre-specified integration accuracy 
by choosing smaller step sizes~\citep{betts_mesh_1998,binder_dynamic_2000,schlegel_dynamic_2005}, 
selecting higher-order local interpolants~\citep{babuska_p-_1990,darby_hp-adaptive_2011},
or relocating grid points~\citep{tanartkit_nested_1997,zhao_density_2011}.
As refinement criteria, 
besides 
minimizing estimates of the integration error,
the adjoint sensitivities can also be used
to guide the refinement process towards lower open-loop costs~\citep{tanartkit_nested_1997,paiva_adaptive_2015}.
The former 
has
been applied to MPC by, e.g., \cite{gao_accurate_2023}, 
the latter, by \cite{lee_mesh_2018} and \cite{paiva_optimal_2018}. 
Tailored to the receding-horizon strategy, 
\cite{paiva_sampleddata_2017} and \cite{potena_non-linear_2018}
propose to save computation time
by enforcing a
higher discretization accuracy 
only for the initial part of the time horizon. 
Nevertheless, due to the variable number of integration steps, the applicability of mesh refinement methods for online applications with a fixed sampling time can be complicated by their unpredictable solve times.
Alternative to rendering the discretization error negligibly small by a suitable discretization strategy, 
\cite{kogel_discrete-time_2015} propose to robustify the prediction against the worst-case integration error by employing a robust \ac{MPC} strategy based on a Lipschitz bound of the dynamics subject to all possible input values; however, this might lead to unnecessarily conservative predictions.

Similarly to \cite{kogel_discrete-time_2015}, rather than accepting a possibly large integration error depending on the optimal input trajectory, or ensuring a specific error tolerance at the cost of unpredictable solve times, in this work we modify the input in response to non-negligible integration errors. 
Instead of a robust treatment, we propose to utilize probabilistic \ac{ODE} solvers to predict the integration error as a function of the applied input. Encoding the trade-off between accurately solving the \ac{ODE} and minimizing the control objective in a stochastic \ac{OCP} formulation, uncertainty reduction can thus be actively targeted by the solver insofar as it is beneficial for the control objective.
Viewing the proposed approach through the lens of probabilistic numerics thereby naturally exposes the implicit assumption of zero discretization error for classical direct transcription methods, as well as provides a dual-control interpretation of the proposed approach in terms of the informativeness of the numerical data processed by the probabilistic \ac{ODE} solver.

To this end, we revisit the continuous-time \ac{OCP} formulation and give an introduction to probabilistic \ac{ODE} solvers in Section~\ref{sec:ProblemFormulation}. In Section~\ref{sec:NewOCP}, we modify the \ac{OCP} to allow for incorporating the integration error quantification provided by the probabilistic \ac{ODE} solver, 
discussing the potential as well as limitations of the proposed approach.
In Section~\ref{sec:NumericalExample}, we illustrate
the proposed method with a numerical example.

\section{Problem formulation} \label{sec:ProblemFormulation}

We consider the continuous-time \ac{OCP}
\begin{subequations} \label{eq:OCP_continuous}
  \begin{align}
    \min_{u(\cdot), x(\cdot)} \quad & \int_0^T \phi(x(t),u(t)) \mathrm{d}t\label{eq:OCP_continuous_cost} \\%  + \phi_f(x(T))  %
    \mathrm{s.t.} \quad & x(0) = x_0 \label{eq:OCP_continuous_initial} \\
    & \dot{x}(t) = f(t, x(t), u(t)), \forall t \in [0,T] \label{eq:OCP_continuous_dynamics} %
  \end{align}
\end{subequations}
where the aim is to find an open-loop input trajectory $u(t)$ minimizing the cost functional~\eqref{eq:OCP_continuous_cost}, subject to nonlinear dynamics given by the \ac{IVP}~\eqref{eq:OCP_continuous_initial},~\eqref{eq:OCP_continuous_dynamics}. 
The solution of the IVP is 
assumed to be
unique and $(p+1)$-times continuously differentiable; 
this is achieved, e.g., when the vector field~$f$ is globally Lipschitz continuous and when $f$ and $u$ are 
$p$-times continuously differentiable, see, e.g., \cite[Thm.~36.6]{hennig2022probabilistic}. 
The stage cost is 
considered
quadratic and positive definite, i.e., \mbox{$\phi(x(t), u(t)) = \frac{1}{2} \left( x(t)^\top W_x x(t) + u(t)^\top W_u u(t) \right)$}, with positive definite matrices $W_x \in \mathbb{R}^{n_x \times n_x}$ and $W_u \in \mathbb{R}^{n_u \times n_u}$.
\begin{remark}
  A quadratic and positive definite terminal cost \mbox{$\phi_f(x(T)) = x(T)^\top P x(T)$} has been omitted from the cost~\eqref{eq:OCP_continuous_cost} in the following exposition without loss of generality.
\end{remark}

\begin{remark}
  In principle, our proposed method also allows for path constraints $h(x(t), u(t)) \leq 0$ to be considered in the optimal control problem. Approaches for including constraints in the \ac{OCP} formulation are discussed in Section~\ref{sec:NewOCP}.
\end{remark}

As the \ac{OCP}~\eqref{eq:OCP_continuous} generally does not admit an analytical solution, we employ the direct transcription method to obtain a numerical approximation. 
For ease of exposition, we choose the single-shooting formulation of the \ac{OCP} and eliminate the state variables. 
Therefore, the input function is parameterized by a parameter vector $\theta := (\theta_1, \ldots, \theta_{N-1})$ in terms of a time-dependent input policy $u(t) = \pi_\theta(t)$ defined on the subintervals $t \in [t_i, t_{i+1})$ of the grid 
$\{ t_i \}_{i=0}^N \subset [0, T]$,
with $0 = t_0 < \ldots < t_N = T$.
By splitting the cost integral into the same subintervals
as well as enforcing the path constraints only 
on the same grid,
we obtain
the direct single-shooting discretization
of~\eqref{eq:OCP_continuous},
  \begin{align} \label{eq:OCP_continuous_q}
      \min_{\theta} \quad & \sum_{i=0}^{N-1} \Phi_i(x_\theta(\cdot),\pi_\theta(\cdot)) %
  \end{align}
where $x_\theta(t_i)$
denotes the \textit{exact} solution to the continuous-time \ac{IVP}~\eqref{eq:OCP_continuous_initial},~\eqref{eq:OCP_continuous_dynamics}
given the control input function \mbox{$u(t) = \pi_\theta(t)$},
and the stage cost $\Phi: \mathbb{R}^{n_x} \times \mathbb{R}^{n_u} \rightarrow \mathbb{R}^{}$ is given as
\begin{align} \label{eq:OCP_continuous_cost_discretization}
  \Phi_i(x(\cdot),u(\cdot)) = \int_{t_i}^{t_{i+1}} \phi(x(t),u(t)) \mathrm{d}t.
\end{align}

Lastly, 
the exact solution of the \ac{IVP}~\eqref{eq:OCP_continuous_initial},~\eqref{eq:OCP_continuous_dynamics}
as well as computation of the integral~\eqref{eq:OCP_continuous_cost_discretization} 
are approximated
in order to obtain a tractable \ac{OCP} formulation.
In this paper, the cost integral is approximated by the Riemann discretization of the quadratic stage cost, that is 
\begin{align}
  \Phi_i(x(\cdot), u(\cdot)) \approx \hat{\Phi}_i(x(t_i), u(t_i)) = \frac{\Delta t_i}{2} \left( x(t_i)^\top W_x x(t_i) + u(t_i)^\top W_u u(t_i) \right), \label{eq:OCP_continuous_cost_quadrature}
\end{align}
where $\Delta t_i \doteq t_{i+1} - t_i$.
For the \ac{ODE} solution, in the following we investigate classical and probabilistic approximation methods through the lens of probabilistic numerics~\citep{hennig2022probabilistic}.

\subsection*{Classical \ac{ODE} solvers}

The solution of the \ac{IVP}~\eqref{eq:OCP_continuous_initial},~\eqref{eq:OCP_continuous_dynamics} as part of the optimization problem is typically performed using classical numerical integration methods, such as Runge--Kutta or collocation methods~\citep{hairer_solving_1993}. For the most part, these methods construct a point estimate
\begin{align}
    \hat{x}_\theta(t_i) \approx x_\theta(t_i)
\end{align}
by extrapolating a $p$-th order polynomial based on the initial condition at the last time step of the approximated solution,
as well as $p \geq 1$ %
evaluations
of the \ac{ODE} vector field. %
Crucially, the extrapolation is performed under the assumption of a perfectly known initial condition, complicating the propagation of integration error across multiple integration steps~\citep{schober_probabilistic_2014}.

In this paper, we employ probabilistic \ac{ODE} solvers to obtain a numerical approximation of the \ac{ODE} solution as part of the \ac{OCP}~\eqref{eq:OCP_continuous}. To this end, in the following section we present a brief introduction to probabilistic \ac{ODE} solvers.

\subsection*{Probabilistic \ac{ODE} solvers}

Probabilistic numerics refers to a class of methods that attempt
to quantify errors in numerical computations probabilistically \citep{hennig_probabilistic_2015,Oates2019a,hennig2022probabilistic}.
Thus, the output of a probabilistic solver to a numerical problem is a probability distribution over candidate solutions rather than a single point.
Probabilistic solvers can roughly be put into two categories:
methods that represent the probability distribution over candidate solutions implicitly via stochasticity \citep{Teymur2016,Teymur2018a,Conrad2017,Abdulle2020},
and methods that represent the distribution over solution candidates explicitly, via the Bayesian formulation \citep{cockayne_bayesian_2019,schober_probabilistic_2014,Kersting2016,Schober2019,tronarp_probabilistic_2019,tronarp_bayesian_2021,kersting_uncertainty-aware_2021}.
This classification is not rigid: Indeed, there are methods that draw on ideas from both concepts \citep{Chkrebtii2016}.

In the Bayesian framework, a solver is constructed by first defining a prior distribution on the quantity of interest, $x(t)$,
and then defining a likelihood reflecting the constraints on $x(t)$ implied by the numerical problem.
In particular, for Bayesian \ac{ODE} solvers, the likelihood is usually defined by subsampling the dynamics constraint~\eqref{eq:OCP_continuous_dynamics},
though this can be extended to other information operators, such as sensor data, algebraic constraints and higher-order derivative information~\citep{bosch_pick-and-mix_2022}.
The aim is thus, given some prior on $x(t)$, to compute a posterior distribution of the form
\begin{align}
  p \left( x(t) \mid x(0) = x_0, \{\dot{x}(t_i) = f(t_i, x(t_i), u(t_i)) \}_{i=0}^N \right),
\end{align}
where, for notational simplicity, 
we present the formulas using the same time discretization for the numerical \ac{ODE} solver and the cost discretization.
In this paper we consider the Gaussian state estimation approach to probabilistic \ac{ODE} solvers, as introduced by
\citet{kersting_active_2016,tronarp_probabilistic_2019,Schober2019};
see also \citet{hennig_probabilistic_2015} for a thorough introduction.
We briefly describe this approach in the following and formulate the numerical \ac{IVP} solution as a Bayesian state estimation problem.

Bayesian \ac{ODE} solvers jointly model our belief about the \ac{ODE} solution $x(t)$ together with its first $p$ derivatives with an extended state 
\begin{align}
  X(t) = \left( x(t), \dot{x}(t), \ldots, \overset{(p)}{x}(t) \right). \label{eq:PN_extended_state}
\end{align}
Restricted to the discretization grid \(\{t_i\}_{i=0}^N\), the 
prior  
is fully described by the density
\begin{align}
  p(X(t_{0:N}) \mid \theta, \kappa) &= \prod_{i=0}^{N-1} p\left( X(t_{i+1}) \mid X(t_i), \theta, \kappa \right) p(X(0) \mid \theta, \kappa)
  \intertext{with the initial and transition densities}
  p(X(0) \mid \theta, \kappa) &= \delta(\bar{X}_0), \\
  p\left( X(t_{i+1}) \mid X(t_i), \theta, \kappa \right) &= \mathcal{N} \left( X(t_{i+1}); A(\Delta t_i) X(t_i), \kappa Q(\Delta t_i) \right),
  \label{eq:transition-densities}
\end{align}
respectively, where
$A(\Delta t_i)$ and $Q(\Delta t_i)$ 
are 
given by the exact discretization of the transition and diffusion matrix of the continuous-time integrated Wiener process prior%
~\cite[see, e.g.,][for the formulas]{kersting_convergence_2020},
and \(\kappa\) is a hyper-parameter scaling the prior uncertainty. 
Note that choosing the prior 
based on a given linear system 
can improve the numerical efficiency for semi-linear dynamics~\citep{bosch_probabilistic_2023-1}.
The initial condition of the extended state $\bar{X}_0$ is chosen such that it satisfies the given initial condition of the \ac{IVP}~\eqref{eq:OCP_continuous_initial},~\eqref{eq:OCP_continuous_dynamics}
for all derivative components
exactly,
and can be computed efficiently via Taylor-mode automatic differentiation \citep{kramer_stable_2020}.

The likelihood model 
is given by the differential equation.
More precisely, we condition the prior process on noise-free observations
\mbox{$z_i := \dot{x}(t_i) - f(t_i, x(t_i), u(t_i)) = 0$}, $i=1,\ldots,N$, which can equivalently be described in terms of the 
extended state as
\begin{align}
  Z(t_i) &= E_1 X(t_i) - f(t_i, E_0 X(t_i), \pi_\theta(t_i)) = 0, \label{eq:PN_observations}
\end{align}
where the matrices $E_0$, $E_1 \in \mathbb{R}^{n_x \times (p+1)n_x}$ are selection matrices extracting $x(t)$ and $\dot{x}(t)$, respectively, from $X(t)$.

A probabilistic estimate of the \ac{ODE} solution and its derivatives is then obtained by computing the posterior density $p(X(t_{0:N}) \mid Z(t_{1:N}), \theta, \kappa)$ according to the following filtering problem: 
\begin{subequations} \label{eq:PN_filtering_problem}
  \begin{align}
    p(X(t_0) \mid \theta, \kappa) &= \delta(\bar{X}_0), \\
    p(X(t_i) \mid X(t_{i-1}), \theta, \kappa) &= \mathcal{N}(X(t_i); A(\Delta t_{i-1}) X(t_{i-1}), \kappa Q(\Delta t_{i-1})), \\
    p(Z(t_i) \mid X(t_i), \theta, \kappa) &= \mathcal{N}(Z(t_i); E_1 X(t_i) - f(E_0 X(t_i), \pi_\theta(t_i)), 0),
  \end{align}
\end{subequations}
for $i=1,\ldots,N$. For linear (affine) vector fields, exact inference can be performed using a \acs{RTS} smoother~\citep{rauch_maximum_1965}, but in the general nonlinear case exact inference is intractable.
We therefore resort to efficient linearization-based approximate inference methods such as the \ac{EKF} 
or the \ac{(I)EKS}. 
Finally, we obtain the Gaussian posterior time marginals
\begin{align}
  \gamma_{\text{prob}}(x_i) \doteq \mathcal{N}(x_i; E_0 \xi_\theta(t_i), \kappa E_0 \Lambda_\theta(t_i) E_0^\top) \approx p(E_0 X(t_i) \mid Z(t_{1:N}), \theta, \kappa), \label{eq:gamma_prob}
\end{align} 
where
the formulas 
for the smoothing mean~$\xi_\theta(t_i)$ and covariance~$\Lambda_\theta(t_i)$
obtained by the 
\ac{(I)EKS}
can be found in appendix~\ref{sec:smoothing_equations}; the calibration of $\kappa$, in appendix~\ref{sec:calibration}.

\begin{remark}[Convergence rates and smoothness of control input]
  When the true solution to 
  the IVP~\eqref{eq:OCP_continuous_initial},~\eqref{eq:OCP_continuous_dynamics} is of smoothness $p+1$,
  for a $p$-times integrated Wiener process prior 
  probabilistic \ac{ODE} solvers are expected to enjoy a global convergence rate of $\mathcal{O}(\delta^p)$~\citep{kersting_convergence_2020,tronarp_bayesian_2021}, where $\delta \doteq \max \{ \Delta t_i \}_{i=0}^{N-1}$.
  The smoothness of the solution is limited by the smoothness of the vector field \citep{Arnold1992}, which in turn is limited by the design of the control signal. 
  Thus,
  it might be necessary to 
  select
  a sufficiently smooth 
  control function 
  to guarantee 
  the desired
  convergence 
  order.

\end{remark}

\section{Probabilistic ODE solvers for optimal control problems} \label{sec:NewOCP}

In order to utilize the posterior uncertainty estimate of the approximate \ac{IVP} solution returned by the probabilistic \ac{ODE} solver, the \ac{OCP} formulation needs to be adapted.
We obtain a well-defined, deterministic cost expression by minimizing the expectation of the cost quadrature with respect to the probabilistic belief over the \ac{ODE} flow, i.e., by solving $\min_\theta V(\theta)$, where
\begin{align}
  V(\theta) &\doteq \mathbb{E}_{x_{0:N-1}} \left[ \sum_{i=0}^{N-1} \hat{\Phi}_i(x_i,\pi_\theta(t_i)) \right], \\ %
  &= \sum_{i=0}^{N-1} \mathbb{E}_{x_i} \left[ \hat{\Phi}_i(x_i,\pi_\theta(t_i)) \right], \\
  &= \sum_{i=0}^{N-1} \int \hat{\Phi}_i(x_i,\pi_\theta(t_i)) \gamma(x_i) \mathrm{d}x_i. \label{eq:expected_cost_integral}
\end{align}
Using the analytical expectation of the quadratic cost in~\eqref{eq:OCP_continuous_cost_quadrature}, we propose to approximate the solution to~\eqref{eq:OCP_continuous_q} with respect to the numerical approximation of the \ac{IVP} solution by solving the \ac{OCP}
\begin{align} \label{eq:OCP_expected}
  \min_{\theta} \quad & \sum_{i=0}^{N-1} \frac{\Delta t_i}{2} \left( \| E_0 \xi_{\theta} (t_i) \|^2_{W_x}  + \| \pi_\theta(t_i) \|^2_{W_u} + \hat{\kappa}_\theta \cdot \mathrm{tr} \left\{W_x E_0 \Lambda_\theta(t_i) E_0^\top \right\} \right), %
\end{align}
where $\xi_{\theta}(t_i)$ and $\Lambda_\theta(t_i)$ are the smoothing mean and covariance of the \acs{EKS}, and $\hat{\kappa}_\theta$ is the prior diffusion scale calibrated according to Appendix~\ref{sec:calibration}, as a function of the inputs~$\theta$.

\begin{remark}
  For a more accurate solution of the filtering problem~\eqref{eq:PN_filtering_problem} the \acs{IEKS} estimate of the filtering problem~\eqref{eq:PN_filtering_problem} can be used by including the \acs{IEKS} fixed-point equations into~\eqref{eq:OCP_expected},
  see appendix~\ref{sec:smoothing_equations}. 
\end{remark}

The reformulation of the deterministic cost in~\eqref{eq:OCP_continuous_q} as an expectation over the probabilistic belief of the \ac{ODE} solution lends itself to the following interpretation from the viewpoint of probabilistic numerics. Clearly, by assigning a Dirac measure to the exact \ac{IVP} solution at time~$t_i$,
\begin{align}
  \gamma_{\text{exact}}(x_i) &\doteq \delta (x_i - x_\theta(t_i)),
  \intertext{minimizing the expected cost in equation~\eqref{eq:expected_cost_integral} with $\gamma(x_i)$ set to $\gamma_{\text{exact}}(x_i)$ still recovers the exact solution to~\eqref{eq:OCP_continuous_q}. Yet, when computing an approximation to the \ac{IVP} solution, we can use the same expression, but instead replace the exact state with the approximation computed by a classical method or a probabilistic \ac{ODE} solver. For classical \ac{ODE} solvers, it becomes evident that likewise assigning the Dirac density shifted by the computed point estimate $\hat{x}_\theta(t_i) \approx x_\theta(t_i)$, i.e.,
  }
  \gamma_{\text{clas}}(x_i) &\doteq \delta (x_i - \hat{x}_\theta(t_i)),
\end{align}
leads to a classical direct single-shooting formulation of the continuous-time \ac{OCP}, where the point estimate~$\hat{x}_\theta(t_i)$ for the \ac{IVP} is simply inserted into the cost~\eqref{eq:OCP_continuous_cost_quadrature}. 
This formalizes the 
standard assumption of zero discretization error employed by classical direct transcription methods.
By including the (approximate)
posterior uncertainty description~\eqref{eq:gamma_prob} provided by the probabilistic \ac{ODE} solver, i.e., propagating the computational uncertainty from the approximate \ac{IVP} solution in the cost of the \ac{OCP}~\eqref{eq:OCP_continuous_q},
this assumption can thus be (partly) relaxed.

Since the uncertainty estimate of the probabilistic \ac{ODE} solver depends on the input trajectory, 
an optimal solution to~\eqref{eq:OCP_expected} captures the trade-off 
between obtaining an optimal mean prediction and reducing the computational uncertainty associated with it. 
Reduction of the computational uncertainty can be achieved by choosing a state-input trajectory that leads to a bigger uncertainty reduction 
in terms of
the posterior covariance $\Lambda_\theta(t_i)$ in~\eqref{eq:OCP_expected}. 
In that sense, the computed optimal inputs 
can exhibit a dual-control effect~\citep{bar-shalom_dual_1974,mesbah_stochastic_2018}, in particular, for nonlinear ODEs.
The connection of the proposed approach to dual control becomes clear when viewed through the lens of probabilistic numerics:
Interpreting the \ac{ODE} constraint~\eqref{eq:PN_observations} as (numerical) data observed within the prediction horizon, the proposed approach can be 
seen
as a stochastic output-feedback MPC acting on the extended state~$X(t)$ in~\eqref{eq:PN_extended_state} conditioned on the data.
However, compared to the setting of most dual control approaches that incorporate the fact that (physical) data is going to be observed in the future, the (numerical) data in this setting is known \emph{a priori}.

The potential of the proposed method in real-time applications with limited computational budget strongly depends on the computational overhead added by solving~\eqref{eq:OCP_expected} with a probabilistic \ac{ODE} solver instead of solving~\eqref{eq:OCP_continuous_q} with a classical \ac{ODE} solver.
Concerning the choice of \ac{ODE} solver, while both classical and probabilistic \ac{ODE} solvers share the same linear time complexity, the complexity in the state dimension for a vanilla probabilistic \ac{ODE} solver implementation is cubic, compared to the linear complexity for most classical solvers. By assuming independence between different state dimensions (which classical solvers generally do since correlations are not modeled), probabilistic \ac{ODE} solvers can enjoy linear space complexity as well~\citep{kramer_probabilistic_2022}. 
Still, to which extent the independence assumption and the implied linearization strategy~\citep[Sec.~2.2]{kramer_probabilistic_2022} affects the dynamics-dependent uncertainty reduction targeted in this paper requires further investigation.
Concerning the choice of \ac{OCP} formulation, optimizing over covariance matrices generally incurs a significant computational overhead, either due to increased nonlinearity or increased state dimension of the condensed or non-condensed \ac{OCP} formulation, respectively.  
Going forward, computational efficiency of the proposed \ac{OCP} formulation will thus have to be addressed in order
for
the uncertainty-aware approach 
to be
competitive with the uncertainty reduction achieved by a (fixed) finer time discretization of a ``certainty-equivalent'' method given the same computational budget.
This paper focuses on highlighting the 
potential 
of using probabilistic integrators 
for integration error-aware control, 
motivating the development of 
more efficient
implementations in future work.

\begin{remark}[Path and terminal constraints] \label{remark:constraints}
  Due to the assumed Gaussian density distribution returned by linearization-based solvers such as the \ac{(I)EKS}, discretized path-constraints and terminal constraints $h(x(t_i), u(t_i)) \leq 0$, $i = 0, \ldots, N$ of the original open-loop problem have to be reformulated probabilistically, 
  based on the probabilistic belief about the \ac{IVP} solution. Computationally efficient approaches include linearization-based constraint tightenings on the posterior mean~\citep{van_hessem_stochastic_2006,kouvaritakis_explicit_2010,cannon_stochastic_2011}, or penalties on the expected constraint violation~\citep{messerer_dual-control_2023}; see also~\cite{mesbah_stochastic_2016} for an overview.
    \end{remark}

\section{Numerical example} \label{sec:NumericalExample}

To illustrate the properties of the algorithm, we apply it to a numerical example%
\footnote{The source code for the numerical experiments is available at \doi{10.3929/ethz-b-000673668}.}%
.
For the probabilistic \ac{ODE} solver, we use the JAX~\citep{jax2018github} implementation of a probabilistic \ac{ODE} solver based on an \acs{EKS} provided by the package \texttt{probdiffeq}~\citep{kramer_probdiffeq_2023}. 
The \ac{OCP}~\eqref{eq:OCP_expected} is implemented in JAX in condensed form using a low-rank variant of the BFGS~\citep{nocedal_numerical_2006} optimizer with bounds (L-BFGSB) optimizer, as implemented in the \texttt{jaxopt} package~\citep{blondel_efficient_2022-1}. 
To test the incorporation of the uncertainty in the \ac{IVP} solution separately from the 
integrator choice, 
the proposed approach 
is compared to a point estimate obtained as the mean prediction of the same probabilistic integrator, which is equivalent to solving~\eqref{eq:OCP_expected} without considering the computational uncertainty, i.e., setting $\Lambda_\theta(t_i) = 0$ in~\eqref{eq:OCP_expected}.
Hence, we denote as ``proposed'' and ``classical'' the proposed approach and proxy for the classical approach, respectively.
The predicted state trajectories of 
both approaches
are compared to
a ground-truth simulation
with the respective optimal control input,
obtained by 
simulating the ODE
using
the IDAS~\citep{cao_adjoint_2003,hindmarsh_sundials_2005} adaptive integration method with a relative error tolerance set to~$10^{-10}$. 

We consider the logistic \ac{ODE} with an added input, 
\mbox{$\dot{x}_1(t) = (1 - x_2(t))^2 x_1(t) - x_2(t) + u(t)$},
\mbox{$\dot{x}_2(t) = x_1(t)$},
with an initial condition $x(0) = (3,1)$ and cost matrices $W_x = 50 I_2$, $W_u = 1$. 
The problem is discretized into $N=20$ control intervals for a time horizon of $T = 5$.
Using a piece-wise affine parameterization of the input function, the probabilistic integrator is applied with integration order $p = 1$.

\begin{figure}
  \begin{minipage}[t]{0.495\textwidth}
    \centering
      \includegraphics{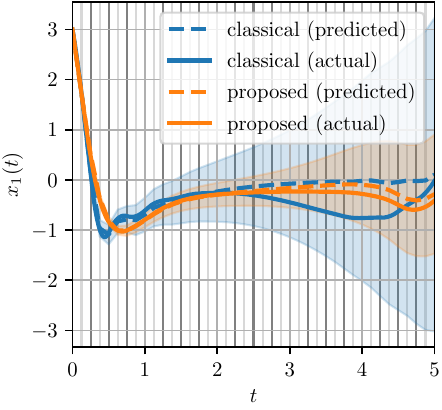}
  \end{minipage}
  \hfill
  \begin{minipage}[b]{0.495\textwidth}
    \centering
    \includegraphics{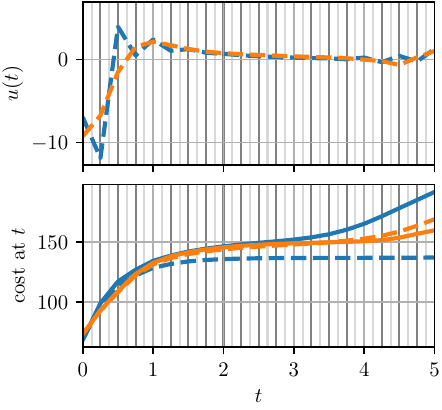}
  \end{minipage}
  \caption{
    \ac{OCP} solutions (dashed lines) compared to ground-truth simulations (solid lines) for the classical (blue) and proposed (orange) approach for $N_\text{int} = 40$
    (left: state $x_1(t)$, top right: cost value, bottom right: applied input). 
    The input computed by the proposed approach 
    actively
    reduces 
    computational uncertainty~(shaded) for a lower expected cost in~\eqref{eq:OCP_expected}, which leads to a more accurate cost prediction and a lower ground-truth cost~\eqref{eq:OCP_continuous_cost_quadrature}.
        }
    \label{fig:XUCost_Plot}
\end{figure}

The simulation results for $N_\text{int} = 40$ integration steps is shown in Figure~\ref{fig:XUCost_Plot}. 
It can be seen that the input computed by the proposed approach (top right image) 
reduces the uncertainty estimate,
leading to a more precise cost estimate (left image), which, in this case, also leads to a lower total actual cost (bottom right image). 
Note that, 
to show the uncertainty reduction achieved by the proposed approach,
the  $2\sigma$-confidence interval of the integrators' marginal uncertainty estimate is also plotted for the classical approach, even though it only considers the mean estimate.

\begin{figure}
  \centering
  \includegraphics{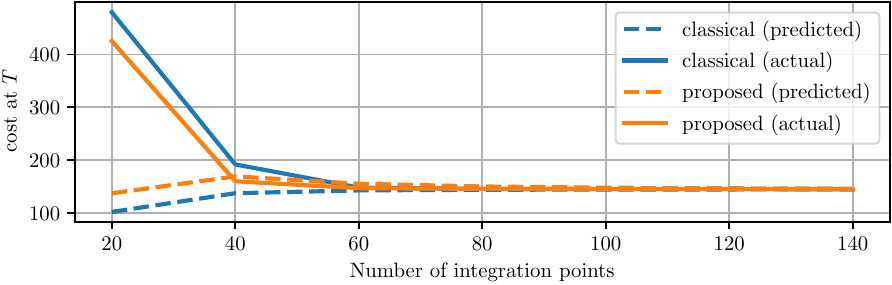}
  \caption{Predicted and ground-truth cost for an increasing number of integration points. 
  For a low number of integration points, the proposed approach achieves a lower actual cost by adapting the input for a more accurate numerical discretization.
  }
  \label{fig:cost_vs_N}
\end{figure}

Figure~\ref{fig:cost_vs_N} shows the predicted and ground-truth cost values at time~$T$ for an increased number of integration points.
Thereby, to isolate the effects of the \ac{ODE} discretization from the cost discretization,
the number of control intervals is kept constant while the number of total integration steps within each control interval is varied.
In the range 
$N_\text{int} \in \{20,40\}$, 
the proposed approach 
leads to a lower actual cost, which is also in better agreement with the predicted cost; for 
$N_\text{int} \geq 60$, 
the slightly underconfident uncertainty estimate leads to a marginally higher predicted cost for the proposed approach, albeit similar actual costs,
while both approaches converge 
as the numerical accuracy increases.
While the proposed approach achieves a better agreement of predicted and actual cost, as well as a lower actual cost for the coarser grids,  
for $N_{\text{int}} = 20$ it can be seen that there is still a considerable difference between the predicted and actual cost value. 
Closer investigations indicate that this discrepancy may be attributed to an overconfident uncertainty calibration; 
the development of improved uncertainty quantification strategies is reserved for future work.

\section{Conclusion and Outlook}

In this paper, we have investigated the 
idea of probabilistic numerics
to
propagate
the numerical uncertainty, 
incurred 
by the approximate discretization of the \ac{IVP} to predict the dynamics inside the \ac{OCP}, 
in the cost prediction.
Quantifying the input-dependent integration error using probabilistic \ac{ODE} solvers,
the computed input can achieve lower ground-truth cost in the presence of computational uncertainty as it automatically determines a trade-off between solving the \ac{IVP} accurately and obtaining a lower nominal cost prediction with respect to the approximate \ac{IVP} solution,
which has been demonstrated in a numerical example. 
For future work, a computationally efficient implementation is to be addressed to achieve these performance benefits in real-time applications with limited computational budget.

\clearpage
\appendix

\section{Technical details on probabilistic \ac{ODE} solvers}

\subsection{Extended Kalman Smoother equations} \label{sec:smoothing_equations}

In the following, we present the equations used for approximating the solution to the smoothing problem~\eqref{eq:PN_filtering_problem}. Thereby we mostly omit the subscript $\theta$ to avoid excessive notation. The \ac{EKF} filtering pass starts from $\mu(t_0) = \bar{X}_0$, $\Sigma(t_0) = 0$ and iterates forward in time over the equations
\begin{subequations}\label{eq:kf}
\begin{align}
\mu(t_i^-)           &= A(\Delta t_{i-1})  \mu(t_{i-1}), \\
\kappa \Sigma(t_i^-) &= \kappa \big( A(\Delta t_{i-1}) \Sigma(t_{i-1}) A^\top(\Delta t_{i-1})  + Q(\Delta t_{i-1})  \big), \\
C(t_i)               &= E_1^\top - E_0^\top J_f(t_i, \tilde{x}(t_i), \pi_\theta(t_i)), \\
b(t_i)               &= f(t_i, \tilde{x}(t_i), \pi_\theta(t_i)) - J_f(t_i, \tilde{x}(t_i), \pi_\theta(t_i)) \tilde{x}(t_i) \\
\kappa S(t_i)        &= \kappa C(t_i) \Sigma(t_i^-) C^\top(t_i), \\
K(t_i)               &= \Sigma(t_i^-) C^\top(t_i) S^{-1}(t_i), \\
\mu(t_i)             &= \mu(t_i^-) +  K(t_i) \big( b(t_i) - C(t_i) \mu(t_i^-)   \big), \\
\kappa \Sigma(t_i)   &= \kappa \big( \Sigma(t_i^-) - K(t_i)S(t_i)K^\top(t_i) \big),
\intertext{
where $J_f$ denotes the Jacobian of $f$ with respect to $x$. The \ac{RTS} smoothing pass starts from the filtering mean and covariance, \mbox{$\xi(t_N) = \mu(t_N)$}, \mbox{$\Lambda(t_N) = \Sigma(t_N)$}, and iterates backwards as follows:%
}
        G(t_i) &= \Sigma(t_i) A^\top(\Delta t_{i+1}) \Sigma^{-1}(t_{i+1}^-), \\
        \kappa P(t_i) &= \kappa(\Sigma(t_i) - G(t_i) \Sigma(t_{i+1}^-) G^\top(t_i)),\\
        \xi(t_i) &= \mu(t_i) + G(t_i) (\xi(t_{i+1}) - \mu(t_{i+1}^-)), \\
        \kappa \Lambda(t_i) &= \kappa (G(t_i) \Lambda(t_{i+1}) G^\top(t_i) + P(t_i)). 
    \end{align}
\end{subequations}
For the \acs{EKS}, the linearization point $\tilde{x}(t_i) = E_0 \mu(t_i^-)$ is the predictive mean of the \ac{EKF}; for the \acs{IEKS}, the smoothing mean $\tilde{x}(t_i) = E_0 \xi(t_i)$ of the previous iteration. 
Note that covariances are invariant to scaling of $\kappa$ \citep{tronarp2019student,tronarp_probabilistic_2019}.
Therefore, the filter can be run for $\kappa = 1$ and the resulting covariances can be post-multiplied once an appropriate value has been determined.

\subsection{Calibration}
\label{sec:calibration}

To return meaningful uncertainty estimates, 
the scaling parameter $\kappa$ needs to be estimated.
We follow the approach described by
\cite{tronarp_probabilistic_2019} and \cite{Bosch2021}
and choose \(\hat{\kappa}\) as the quasi-\ac{MLE}.
In view of \eqref{eq:kf}, this estimate can be obtained by solving
\begin{align}
  \hat{\kappa}_\theta &= \underset{\kappa}{\arg \max} \> \sum_{i=1}^N \log \mathcal{N}(b(t_i); C(t_i) \mu(t_i^-), \kappa S(t_i)) \\
  &= \frac{1}{Nn_x} \sum_{i = 1}^N \big( b(t_i) - C(t_i) \mu(t_i^-)   \big)^\top S^{-1}(t_i) \big( b(t_i) - C(t_i) \mu(t_i^-)   \big).
\end{align}

\acks{
The authors thank Nicholas Krämer for his support in using the \texttt{probdiffeq} software library,
Andrea Carron, Johannes Köhler and Anna Scampicchio for providing helpful feedback,
as well as Moritz Diehl and Armin Nurkanović for inspiring discussions on the initial motivation of the project.
AL and MZ received funding by the European Union's Horizon 2020 research and innovation programme, Marie~Sk\l{}odowska-Curie grant agreement No. 953348,~\mbox{ELO-X}. FT was partially supported by the Wallenberg AI, Autonomous Systems and
Software Program (WASP) funded by the Knut and Alice Wallenberg Foundation.
PH, NB and JS gratefully acknowledge financial support by the DFG Cluster of Excellence "Machine Learning - New Perspectives for Science", EXC 2064/1, project number 390727645; the German Federal Ministry of Education and Research (BMBF) through the Tübingen AI Center (FKZ: 01IS18039A), as well as funds from the Ministry of Science, Research and Arts of the State of Baden-Württemberg. 
The authors also thank the International Max Planck Research School for Intelligent Systems (IMPRS-IS) for supporting NB and JS.
}

\bibliography{bibliography}

\begin{thebibliography}{63}
\providecommand{\natexlab}[1]{#1}
\providecommand{\url}[1]{\texttt{#1}}
\expandafter\ifx\csname urlstyle\endcsname\relax
  \providecommand{\doi}[1]{doi: #1}\else
  \providecommand{\doi}{doi: \begingroup \urlstyle{rm}\Url}\fi

\bibitem[Abdulle and Garegnani(2020)]{Abdulle2020}
A.~Abdulle and G.~Garegnani.
\newblock Random time step probabilistic methods for uncertainty quantification
  in chaotic and geometric numerical integration.
\newblock \emph{Statistics and Computing}, 2020.

\bibitem[Arnold(1992)]{Arnold1992}
Vladimir~I Arnold.
\newblock \emph{Ordinary Differential Equations}.
\newblock Springer Science \& Business Media, 1992.

\bibitem[Babu{\v s}ka and Suri(1990)]{babuska_p-_1990}
Ivo Babu{\v s}ka and Manil Suri.
\newblock The p- and h-p versions of the finite element method, an overview.
\newblock \emph{Computer Methods in Applied Mechanics and Engineering},
  80\penalty0 (1), 1990.
\newblock \doi{10.1016/0045-7825(90)90011-A}.

\bibitem[{Bar-Shalom} and Tse(1974)]{bar-shalom_dual_1974}
Y.~{Bar-Shalom} and E.~Tse.
\newblock Dual effect, certainty equivalence, and separation in stochastic
  control.
\newblock \emph{IEEE Transactions on Automatic Control}, 19\penalty0 (5), 1974.
\newblock \doi{10.1109/TAC.1974.1100635}.

\bibitem[Betts and Huffman(1998)]{betts_mesh_1998}
John~T. Betts and William~P. Huffman.
\newblock Mesh refinement in direct transcription methods for optimal control.
\newblock \emph{Optimal Control Applications and Methods}, 19\penalty0 (1),
  1998.
\newblock
  \doi{10.1002/(SICI)1099-1514(199801/02)19:1<1::AID-OCA616>3.0.CO;2-Q}.

\bibitem[Binder et~al.(2000)Binder, Cruse, Cruz~Villar, and
  Marquardt]{binder_dynamic_2000}
T.~Binder, A.~Cruse, C.A. Cruz~Villar, and W.~Marquardt.
\newblock Dynamic optimization using a wavelet based adaptive control vector
  parameterization strategy.
\newblock \emph{Computers \& Chemical Engineering}, 24\penalty0 (2-7), 2000.
\newblock \doi{10.1016/S0098-1354(00)00357-4}.

\bibitem[Blondel et~al.(2022)Blondel, Berthet, Cuturi, Frostig, Hoyer,
  {Llinares-L{\'o}pez}, Pedregosa, and Vert]{blondel_efficient_2022-1}
Mathieu Blondel, Quentin Berthet, Marco Cuturi, Roy Frostig, Stephan Hoyer,
  Felipe {Llinares-L{\'o}pez}, Fabian Pedregosa, and Jean-Philippe Vert.
\newblock Efficient and {{Modular Implicit Differentiation}}, 2022.

\bibitem[Bosch et~al.(2021)Bosch, Hennig, and Tronarp]{Bosch2021}
N.~Bosch, P.~Hennig, and F.~Tronarp.
\newblock Calibrated adaptive probabilistic {{ODE}} solvers.
\newblock In \emph{24th International Conference on Artificial Intelligence and
  Statistics}, 2021.

\bibitem[Bosch et~al.(2022)Bosch, Tronarp, and Hennig]{bosch_pick-and-mix_2022}
Nathanael Bosch, Filip Tronarp, and Philipp Hennig.
\newblock Pick-and-{{Mix Information Operators}} for {{Probabilistic ODE
  Solvers}}.
\newblock In \emph{Proceedings of {{The}} 25th {{International Conference}} on
  {{Artificial Intelligence}} and {{Statistics}}}. PMLR, 2022.
\newblock URL \url{https://proceedings.mlr.press/v151/bosch22a.html}.

\bibitem[Bosch et~al.(2023)Bosch, Hennig, and
  Tronarp]{bosch_probabilistic_2023-1}
Nathanael Bosch, Philipp Hennig, and Filip Tronarp.
\newblock Probabilistic {{Exponential Integrators}}.
\newblock \emph{Advances in Neural Information Processing Systems}, 36, 2023.
\newblock URL
  \url{https://proceedings.neurips.cc/paper_files/paper/2023/hash/7f64034009f4a5fa417a57e1a987c5cd-Abstract-Conference.html}.

\bibitem[Bradbury et~al.(2018)Bradbury, Frostig, Hawkins, Johnson, Leary,
  Maclaurin, Necula, Paszke, VanderPlas, {Wanderman-Milne}, and
  Zhang]{jax2018github}
James Bradbury, Roy Frostig, Peter Hawkins, Matthew~James Johnson, Chris Leary,
  Dougal Maclaurin, George Necula, Adam Paszke, Jake VanderPlas, Skye
  {Wanderman-Milne}, and Qiao Zhang.
\newblock {{JAX}}: Composable transformations of {{Python}}+{{NumPy}} programs,
  2018.
\newblock URL \url{http://github.com/google/jax}.

\bibitem[Cagienard et~al.(2007)Cagienard, Grieder, Kerrigan, and
  Morari]{cagienard_move_2007}
R.~Cagienard, P.~Grieder, E.~C. Kerrigan, and M.~Morari.
\newblock Move blocking strategies in receding horizon control.
\newblock \emph{Journal of Process Control}, 17\penalty0 (6), 2007.
\newblock \doi{10.1016/j.jprocont.2007.01.001}.

\bibitem[Cannon et~al.(2011)Cannon, Kouvaritakis, Rakovi{\'c}, and
  Cheng]{cannon_stochastic_2011}
Mark Cannon, Basil Kouvaritakis, Sa{\v s}a~V. Rakovi{\'c}, and Qifeng Cheng.
\newblock Stochastic {{Tubes}} in {{Model Predictive Control With Probabilistic
  Constraints}}.
\newblock \emph{IEEE Transactions on Automatic Control}, 56\penalty0 (1), 2011.
\newblock \doi{10.1109/TAC.2010.2086553}.

\bibitem[Cao et~al.(2003)Cao, Li, Petzold, and Serban]{cao_adjoint_2003}
Yang Cao, Shengtai Li, Linda Petzold, and Radu Serban.
\newblock Adjoint {{Sensitivity Analysis}} for {{Differential-Algebraic
  Equations}}: {{The Adjoint DAE System}} and {{Its Numerical Solution}}.
\newblock \emph{SIAM J. Sci. Comput.}, 24\penalty0 (3), 2003.
\newblock \doi{10.1137/S1064827501380630}.

\bibitem[Chen et~al.(2020)Chen, Scarabottolo, Bruschetta, and
  Beghi]{chen_efficient_2020}
Yutao Chen, Nicol{\'o} Scarabottolo, Mattia Bruschetta, and Alessandro Beghi.
\newblock Efficient move blocking strategy for multiple shooting-based
  non-linear model predictive control.
\newblock \emph{IET Control Theory \& Applications}, 14\penalty0 (2), 2020.
\newblock \doi{10.1049/iet-cta.2019.0168}.

\bibitem[Chkrebtii et~al.(2016)Chkrebtii, Campbell, Calderhead, and
  Girolami]{Chkrebtii2016}
O.~A. Chkrebtii, D.~A. Campbell, B.~Calderhead, and M.~A. Girolami.
\newblock Bayesian solution uncertainty quantification for differential
  equations.
\newblock \emph{Bayesian Analysis}, 11\penalty0 (4), 2016.

\bibitem[Cockayne et~al.(2019)Cockayne, Oates, Sullivan, and
  Girolami]{cockayne_bayesian_2019}
Jon Cockayne, Chris~J. Oates, T.~J. Sullivan, and Mark Girolami.
\newblock Bayesian {{Probabilistic Numerical Methods}}.
\newblock \emph{SIAM Rev.}, 61\penalty0 (4), 2019.
\newblock \doi{10.1137/17M1139357}.

\bibitem[Conrad et~al.(2017)Conrad, Girolami, S{\"a}rkk{\"a}, Stuart, and
  Zygalakis]{Conrad2017}
P.~R. Conrad, M.~Girolami, S.~S{\"a}rkk{\"a}, A.~Stuart, and K.~Zygalakis.
\newblock Statistical analysis of differential equations: Introducing
  probability measures on numerical solutions.
\newblock \emph{Statistics and Computing}, 27\penalty0 (4), 2017.

\bibitem[Darby et~al.(2011)Darby, Hager, and Rao]{darby_hp-adaptive_2011}
Christopher~L. Darby, William~W. Hager, and Anil~V. Rao.
\newblock An hp-adaptive pseudospectral method for solving optimal control
  problems.
\newblock \emph{Optimal Control Applications and Methods}, 32\penalty0 (4),
  2011.
\newblock \doi{10.1002/oca.957}.

\bibitem[Gao et~al.(2023)Gao, Han, Eben~Li, Xu, and Dang]{gao_accurate_2023}
Feng Gao, Yu~Han, Shengbo Eben~Li, Shaobing Xu, and Dongfang Dang.
\newblock Accurate {{Pseudospectral Optimization}} of {{Nonlinear Model
  Predictive Control}} for {{High-Performance Motion Planning}}.
\newblock \emph{IEEE Transactions on Intelligent Vehicles}, 8\penalty0 (2),
  2023.
\newblock \doi{10.1109/TIV.2022.3153633}.

\bibitem[Hairer et~al.(1993)Hairer, Wanner, and
  N{\o}rsett]{hairer_solving_1993}
Ernst Hairer, Gerhard Wanner, and Syvert~P. N{\o}rsett.
\newblock \emph{Solving {{Ordinary Differential Equations I}}}, volume~8 of
  \emph{Springer {{Series}} in {{Computational Mathematics}}}.
\newblock Springer, Berlin, Heidelberg, 1993.
\newblock ISBN 978-3-540-56670-0 978-3-540-78862-1.
\newblock \doi{10.1007/978-3-540-78862-1}.

\bibitem[Hennig et~al.(2015)Hennig, Osborne, and
  Girolami]{hennig_probabilistic_2015}
Philipp Hennig, Michael~A. Osborne, and Mark Girolami.
\newblock Probabilistic numerics and uncertainty in computations.
\newblock \emph{Proc. R. Soc. A.}, 471\penalty0 (2179), 2015.
\newblock \doi{10.1098/rspa.2015.0142}.

\bibitem[Hennig et~al.(2022)Hennig, Osborne, and
  Kersting]{hennig2022probabilistic}
Philipp Hennig, Michael~A Osborne, and Hans~P Kersting.
\newblock \emph{Probabilistic Numerics: {{Computation}} as Machine Learning}.
\newblock Cambridge University Press, 2022.

\bibitem[Hindmarsh et~al.(2005)Hindmarsh, Brown, Grant, Lee, Serban, Shumaker,
  and Woodward]{hindmarsh_sundials_2005}
Alan~C. Hindmarsh, Peter~N. Brown, Keith~E. Grant, Steven~L. Lee, Radu Serban,
  Dan~E. Shumaker, and Carol~S. Woodward.
\newblock {{SUNDIALS}}: {{Suite}} of nonlinear and differential/algebraic
  equation solvers.
\newblock \emph{ACM Trans. Math. Softw.}, 31\penalty0 (3), 2005.
\newblock \doi{10.1145/1089014.1089020}.

\bibitem[Kersting and Hennig(2016{\natexlab{a}})]{Kersting2016}
H.~Kersting and P.~Hennig.
\newblock Active uncertainty calibration in {{Bayesian ODE}} solvers.
\newblock In \emph{Uncertainty in Artificial Intelligence ({{UAI}}) 2016},
  2016{\natexlab{a}}.

\bibitem[Kersting(2021)]{kersting_uncertainty-aware_2021}
Hans Kersting.
\newblock \emph{Uncertainty-{{Aware Numerical Solutions}} of {{ODEs}} by
  {{Bayesian Filtering}}}.
\newblock PhD thesis, Universit{\"a}t T{\"u}bingen, 2021.

\bibitem[Kersting and Hennig(2016{\natexlab{b}})]{kersting_active_2016}
Hans Kersting and Philipp Hennig.
\newblock Active {{Uncertainty Calibration}} in {{Bayesian ODE Solvers}}.
\newblock In \emph{{{UAI}} 2016 {{Proceedings}}}, NY, USA, 2016{\natexlab{b}}.
\newblock URL \url{http://www.auai.org/uai2016/proceedings/papers/163.pdf}.

\bibitem[Kersting et~al.(2020)Kersting, Sullivan, and
  Hennig]{kersting_convergence_2020}
Hans Kersting, T.~J. Sullivan, and Philipp Hennig.
\newblock Convergence rates of {{Gaussian ODE}} filters.
\newblock \emph{Stat Comput}, 30\penalty0 (6), 2020.
\newblock \doi{10.1007/s11222-020-09972-4}.

\bibitem[K{\"o}gel and Findeisen(2015)]{kogel_discrete-time_2015}
Markus K{\"o}gel and Rolf Findeisen.
\newblock Discrete-time robust model predictive control for continuous-time
  nonlinear systems.
\newblock In \emph{2015 {{American Control Conference}} ({{ACC}})}, 2015.
\newblock \doi{10.1109/ACC.2015.7170852}.

\bibitem[Kouvaritakis et~al.(2010)Kouvaritakis, Cannon, Rakovi{\'c}, and
  Cheng]{kouvaritakis_explicit_2010}
Basil Kouvaritakis, Mark Cannon, Sa{\v s}a~V. Rakovi{\'c}, and Qifeng Cheng.
\newblock Explicit use of probabilistic distributions in linear predictive
  control.
\newblock \emph{Automatica}, 46\penalty0 (10), 2010.
\newblock \doi{10.1016/j.automatica.2010.06.034}.

\bibitem[Kr{\"a}mer(2023)]{kramer_probdiffeq_2023}
Nicholas Kr{\"a}mer.
\newblock Probdiffeq, 2023.
\newblock URL \url{https://github.com/pnkraemer/probdiffeq}.

\bibitem[Kr{\"a}mer and Hennig(2020)]{kramer_stable_2020}
Nicholas Kr{\"a}mer and Philipp Hennig.
\newblock Stable {{Implementation}} of {{Probabilistic ODE Solvers}}, 2020.

\bibitem[Kr{\"a}mer et~al.(2022)Kr{\"a}mer, Bosch, Schmidt, and
  Hennig]{kramer_probabilistic_2022}
Nicholas Kr{\"a}mer, Nathanael Bosch, Jonathan Schmidt, and Philipp Hennig.
\newblock Probabilistic {{ODE Solutions}} in {{Millions}} of {{Dimensions}}.
\newblock In \emph{Proceedings of the 39th {{International Conference}} on
  {{Machine Learning}}}. PMLR, 2022.
\newblock URL \url{https://proceedings.mlr.press/v162/kramer22b.html}.

\bibitem[Lazutkin et~al.(2018)Lazutkin, Geletu, and Li]{lazutkin_approach_2018}
Evgeny Lazutkin, Abebe Geletu, and Pu~Li.
\newblock An {{Approach}} to {{Determining}} the {{Number}} of {{Time
  Intervals}} for {{Solving Dynamic Optimization Problems}}.
\newblock \emph{Ind. Eng. Chem. Res.}, 57\penalty0 (12), 2018.
\newblock \doi{10.1021/acs.iecr.7b03361}.

\bibitem[Lee et~al.(2018)Lee, Moase, and Manzie]{lee_mesh_2018}
K.~Lee, W.~H. Moase, and C.~Manzie.
\newblock Mesh adaptation in direct collocated nonlinear model predictive
  control.
\newblock \emph{International Journal of Robust and Nonlinear Control},
  28\penalty0 (15), 2018.
\newblock \doi{10.1002/rnc.4235}.

\bibitem[Lishkova et~al.(2022)Lishkova, Cannon, and
  {Ober-Bl{\"o}baum}]{lishkova_multirate_2022}
Yana Lishkova, Mark Cannon, and Sina {Ober-Bl{\"o}baum}.
\newblock A {{Multirate Variational Approach}} to {{Nonlinear MPC}}.
\newblock In \emph{2022 {{European Control Conference}} ({{ECC}})}, 2022.
\newblock \doi{10.23919/ECC55457.2022.9838240}.

\bibitem[Mesbah(2016)]{mesbah_stochastic_2016}
Ali Mesbah.
\newblock Stochastic {{Model Predictive Control}}: {{An Overview}} and
  {{Perspectives}} for {{Future Research}}.
\newblock \emph{IEEE Control Systems Magazine}, 36\penalty0 (6), 2016.
\newblock \doi{10.1109/MCS.2016.2602087}.

\bibitem[Mesbah(2018)]{mesbah_stochastic_2018}
Ali Mesbah.
\newblock Stochastic model predictive control with active uncertainty learning:
  {{A Survey}} on dual control.
\newblock \emph{Annual Reviews in Control}, 45, 2018.
\newblock \doi{10.1016/j.arcontrol.2017.11.001}.

\bibitem[Messerer et~al.(2023)Messerer, Baumg{\"a}rtner, and
  Diehl]{messerer_dual-control_2023}
Florian Messerer, Katrin Baumg{\"a}rtner, and Moritz Diehl.
\newblock A {{Dual-Control Effect Preserving Formulation}} for {{Nonlinear
  Output-Feedback Stochastic Model Predictive Control With Constraints}}.
\newblock \emph{IEEE Control Systems Letters}, 7, 2023.
\newblock \doi{10.1109/LCSYS.2022.3230552}.

\bibitem[Nocedal and Wright(2006)]{nocedal_numerical_2006}
Jorge Nocedal and Stephen~J. Wright.
\newblock \emph{Numerical Optimization}.
\newblock Springer Series in Operation Research and Financial Engineering.
  Springer, New York, second edition edition, 2006.
\newblock ISBN 0-387-30303-0.

\bibitem[Oates and Sullivan(2019)]{Oates2019a}
C.~J. Oates and T.~J. Sullivan.
\newblock A modern retrospective on probabilistic numerics.
\newblock \emph{Statistics and Computing}, 29\penalty0 (6), 2019.

\bibitem[Paiva and Fontes(2015)]{paiva_adaptive_2015}
Lu{\'i}s~Tiago Paiva and Fernando A. C.~C. Fontes.
\newblock Adaptive time--mesh refinement in optimal control problems with state
  constraints.
\newblock \emph{Discrete \& Continuous Dynamical Systems - A}, 35\penalty0 (9),
  2015.
\newblock \doi{10.3934/dcds.2015.35.4553}.

\bibitem[Paiva and Fontes(2017)]{paiva_sampleddata_2017}
Lu{\'i}s~Tiago Paiva and Fernando A. C.~C. Fontes.
\newblock Sampled--{{Data Model Predictive Control Using Adaptive
  Time}}--{{Mesh Refinement Algorithms}}.
\newblock In \emph{{{CONTROLO}} 2016}, Lecture {{Notes}} in {{Electrical
  Engineering}}, Cham, 2017. Springer International Publishing.
\newblock ISBN 978-3-319-43671-5.
\newblock \doi{10.1007/978-3-319-43671-5_13}.

\bibitem[Paiva and Fontes(2018)]{paiva_optimal_2018}
Lu{\'i}s~Tiago Paiva and Fernando A. C.~C. Fontes.
\newblock Optimal {{Control Algorithms}} with {{Adaptive Time-Mesh Refinement}}
  for {{Kite Power Systems}}.
\newblock \emph{Energies}, 11\penalty0 (3), 2018.
\newblock \doi{10.3390/en11030475}.

\bibitem[Potena et~al.(2018)Potena, Della~Corte, Nardi, Grisetti, and
  Pretto]{potena_non-linear_2018}
Ciro Potena, Bartolomeo Della~Corte, Daniele Nardi, Giorgio Grisetti, and
  Alberto Pretto.
\newblock Non-linear model predictive control with adaptive time-mesh
  refinement.
\newblock In \emph{2018 {{IEEE International Conference}} on {{Simulation}},
  {{Modeling}}, and {{Programming}} for {{Autonomous Robots}} ({{SIMPAR}})},
  2018.
\newblock \doi{10.1109/SIMPAR.2018.8376274}.

\bibitem[Qin and Badgwell(1997)]{qin1997overview}
S~Joe Qin and Thomas~A Badgwell.
\newblock An overview of industrial model predictive control technology.
\newblock In \emph{{{AIche}} Symposium Series}, volume~93. New York, NY:
  American Institute of Chemical Engineers, 1971-c2002., 1997.

\bibitem[Quirynen et~al.(2015)Quirynen, Vukov, and
  Diehl]{quirynen_multiple_2015}
Rien Quirynen, Milan Vukov, and Moritz Diehl.
\newblock Multiple {{Shooting}} in a {{Microsecond}}.
\newblock In \emph{Multiple {{Shooting}} and {{Time Domain Decomposition
  Methods}}}, Contributions in {{Mathematical}} and {{Computational Sciences}},
  Cham, 2015. Springer International Publishing.
\newblock ISBN 978-3-319-23321-5.
\newblock \doi{10.1007/978-3-319-23321-5_7}.

\bibitem[Rauch et~al.(1965)Rauch, Tung, and Striebel]{rauch_maximum_1965}
H.~E. Rauch, F.~Tung, and C.~T. Striebel.
\newblock Maximum likelihood estimates of linear dynamic systems.
\newblock \emph{AIAA Journal}, 3\penalty0 (8), 1965.
\newblock \doi{10.2514/3.3166}.

\bibitem[Schlegel et~al.(2005)Schlegel, Stockmann, Binder, and
  Marquardt]{schlegel_dynamic_2005}
Martin Schlegel, Klaus Stockmann, Thomas Binder, and Wolfgang Marquardt.
\newblock Dynamic optimization using adaptive control vector parameterization.
\newblock \emph{Computers \& Chemical Engineering}, 29\penalty0 (8), 2005.
\newblock \doi{10.1016/j.compchemeng.2005.02.036}.

\bibitem[Schober et~al.(2019)Schober, S{\"a}rkk{\"a}, and Hennig]{Schober2019}
M.~Schober, S.~S{\"a}rkk{\"a}, and P.~Hennig.
\newblock A probabilistic model for the numerical solution of initial value
  problems.
\newblock \emph{Statistics and Computing}, 29\penalty0 (1), 2019.

\bibitem[Schober et~al.(2014)Schober, Duvenaud, and
  Hennig]{schober_probabilistic_2014}
Michael Schober, David~K Duvenaud, and Philipp Hennig.
\newblock Probabilistic {{ODE Solvers}} with {{Runge-Kutta Means}}.
\newblock In \emph{Advances in {{Neural Information Processing Systems}}},
  volume~27. Curran Associates, Inc., 2014.
\newblock URL
  \url{https://proceedings.neurips.cc/paper/2014/hash/59b90e1005a220e2ebc542eb9d950b1e-Abstract.html}.

\bibitem[Shekhar and Manzie(2015)]{shekhar_optimal_2015}
Rohan~C. Shekhar and Chris Manzie.
\newblock Optimal move blocking strategies for model predictive control.
\newblock \emph{Automatica}, 61, 2015.
\newblock \doi{10.1016/j.automatica.2015.07.030}.

\bibitem[Tanartkit and Biegler(1997)]{tanartkit_nested_1997}
P.~Tanartkit and L.~T. Biegler.
\newblock A nested, simultaneous approach for dynamic optimization
  problems---{{II}}: The outer problem.
\newblock \emph{Computers \& Chemical Engineering}, 21\penalty0 (12), 1997.
\newblock \doi{10.1016/S0098-1354(97)00014-8}.

\bibitem[Teymur et~al.(2016)Teymur, Zygalakis, and Calderhead]{Teymur2016}
Onur Teymur, Kostas Zygalakis, and Ben Calderhead.
\newblock Probabilistic linear multistep methods.
\newblock In \emph{Advances in Neural Information Processing Systems
  ({{NIPS}})}, 2016.

\bibitem[Teymur et~al.(2018)Teymur, Lie, Sullivan, and Calderhead]{Teymur2018a}
Onur Teymur, Han~Cheng Lie, Tim Sullivan, and Ben Calderhead.
\newblock Implicit probabilistic integrators for {{ODEs}}.
\newblock In \emph{Advances in Neural Information Processing Systems
  ({{NIPS}})}, 2018.

\bibitem[T{\o}ndel and Johansen(2002)]{tondel_complexity_2002}
Petter T{\o}ndel and Tor~A. Johansen.
\newblock Complexity reduction in explicit linear model predictive control.
\newblock \emph{IFAC Proceedings Volumes}, 35\penalty0 (1), 2002.
\newblock \doi{10.3182/20020721-6-ES-1901.00600}.

\bibitem[Tronarp et~al.(2019{\natexlab{a}})Tronarp, Karvonen, and
  S{\"a}rkk{\"a}]{tronarp2019student}
Filip Tronarp, Toni Karvonen, and Simo S{\"a}rkk{\"a}.
\newblock Student's {{ t }}-{{Filters}} for noise scale estimation.
\newblock \emph{IEEE Signal Processing Letters}, 26\penalty0 (2),
  2019{\natexlab{a}}.

\bibitem[Tronarp et~al.(2019{\natexlab{b}})Tronarp, Kersting, S{\"a}rkk{\"a},
  and Hennig]{tronarp_probabilistic_2019}
Filip Tronarp, Hans Kersting, Simo S{\"a}rkk{\"a}, and Philipp Hennig.
\newblock Probabilistic solutions to ordinary differential equations as
  nonlinear {{Bayesian}} filtering: A new perspective.
\newblock \emph{Stat Comput}, 29\penalty0 (6), 2019{\natexlab{b}}.
\newblock \doi{10.1007/s11222-019-09900-1}.

\bibitem[Tronarp et~al.(2021)Tronarp, S{\"a}rkk{\"a}, and
  Hennig]{tronarp_bayesian_2021}
Filip Tronarp, Simo S{\"a}rkk{\"a}, and Philipp Hennig.
\newblock Bayesian {{ODE}} solvers: The maximum a posteriori estimate.
\newblock \emph{Stat Comput}, 31\penalty0 (3), 2021.
\newblock \doi{10.1007/s11222-021-09993-7}.

\bibitem[Van~Hessem and Bosgra(2006)]{van_hessem_stochastic_2006}
Dennis Van~Hessem and Okko Bosgra.
\newblock Stochastic closed-loop model predictive control of continuous
  nonlinear chemical processes.
\newblock \emph{Journal of Process Control}, 16\penalty0 (3), 2006.
\newblock \doi{10.1016/j.jprocont.2005.06.003}.

\bibitem[Wolff et~al.(2022)Wolff, Uchiyama, Burkhardt, and
  Sawodny]{wolff_nonlinear_2022}
Frank Wolff, Naoki Uchiyama, Mark Burkhardt, and Oliver Sawodny.
\newblock Nonlinear {{Model Predictive Control}} with {{Non-Equidistant
  Discretization Time Grids}} for {{Rotary Cranes}}.
\newblock In \emph{2022 13th {{Asian Control Conference}} ({{ASCC}})}, 2022.
\newblock \doi{10.23919/ASCC56756.2022.9828180}.

\bibitem[Yu and Biegler(2016)]{yu_stable_2016}
Mingzhao Yu and Lorenz~T. Biegler.
\newblock A {{Stable}} and {{Robust NMPC Strategy}} with {{Reduced Models}} and
  {{Nonuniform Grids}}.
\newblock \emph{IFAC-PapersOnLine}, 49\penalty0 (7), 2016.
\newblock \doi{10.1016/j.ifacol.2016.07.212}.

\bibitem[Zhao and Tsiotras(2011)]{zhao_density_2011}
Yiming Zhao and Panagiotis Tsiotras.
\newblock Density {{Functions}} for {{Mesh Refinement}} in {{Numerical Optimal
  Control}}.
\newblock \emph{Journal of Guidance, Control, and Dynamics}, 34\penalty0 (1),
  2011.
\newblock \doi{10.2514/1.45852}.

\end{thebibliography}

\end{document}